\documentclass{article}
\usepackage{latexsym} 
\usepackage[all]{xy}
\usepackage{amsmath}
\usepackage{amsfonts}

\newtheorem{thm}{Theorem}[section]
\newtheorem{prop}[thm]{Proposition}
\newtheorem{lem}[thm]{Lemma}
\newtheorem{df}[thm]{Definition}

\begin{document}

\title{On motives for Deligne-Mumford stacks}
\author{B. Toen\footnote{Max Planck Institut fur Mathematik, Vivatsgasse 7, 53111 Bonn, Germany \newline
 \hspace*{5mm} toen@mpim-bonn.mpg.de}}

\maketitle

\begin{abstract}
We define and compare two different definitions of Chow motives for Deligne-Mumford stacks, associated with two different
definitions of Chow rings. The main result we prove is that both categories of motives are equivalent to the usual category of
motives of algebraic varieties, but the motives of a given stack associated with both theories are not isomorphic. 
We will also give some examples of motives associated with some algebraic stacks. 
\end{abstract}

\tableofcontents

\newpage

The construction of Gromov-Witten invariants in algebraic geometry is based on two fundamental objects. The first one is a 
diagram of algebraic stacks
$$\xymatrix{
 & \overline{\mathcal{M}}_{g,n}(V,\beta) \ar[dl]_-{ev} \ar[rd]^-{o} & \\
V^{n} & & \overline{\mathcal{M}}_{g,n}}$$
where $\overline{\mathcal{M}}_{g,n}(V,\beta)$ is the stack of stable maps $f : C \longrightarrow V$, with $f_{*}([C])=\beta$,
$o$ the morphism which forgets $f$ then stabilizes, and $ev$ the evaluation morphism at the marked points on the curve $C$.

The second one is the virtual fundamental class, $\mathcal{I}_{g,n}(V,\beta) \in A_{D}(\overline{\mathcal{M}}_{g,n}(V,\beta))$ (\cite{bf}). 

These two objects combine to give the Gromov-Witten correspondence (\cite{b})
$$I_{g,n}(V,\beta):=(ev,o)_{*}(\mathcal{I}_{g,n}(V,\beta)) \in A_{D}(V^{n}\times \overline{\mathcal{M}}_{g,n}),$$
from which the Gromov-Witten invariants are defined. This shows the motivic nature of Gromov-Witten invariants, and raises the question of the construction of a good theory of motives for Deligne-Mumford stacks. \\

We possess now at least two different ways to define such a theory, corresponding to two different 
definitions of Chow cohomology for Deligne-Mumford stacks. There exist first the $A^{*}$ theories, which are described 
in \cite{eg,g,j,k,v}, and which all coincide with rational coefficients. They satisfy every expected properties of a Chow
cohomology, except the Riemann-Roch theorem, and in particular the Riemann-Roch isomorphism $K_{0}(F)_{\mathbb{Q}} \simeq
A^{*}(F)_{\mathbb{Q}}$. On the other hand, we have the enriched theory $A^{*}_{\chi}$, which was used in \cite{t1,t2} to prove the 
Grothendieck-Riemann-Roch formula, and for which $K_{0}(F)_{\mathbb{Q}}\simeq A^{*}_{\chi}(F)_{\mathbb{Q}}$. \\

In this paper we will show that the two categories of motives associated to these two previous theories are equivalent. This will be shown
using the fact that both of them are equivalent to the usual category of Chow motives for algebraic varieties
(answering a question of Y. Manin and K. Behrend \cite[$8.2$]{bm}). However this does not implies
that the two associated motivic theories are equivalent. Indeed, we will show that the motive of a Deligne-Mumford stack associated
with the theory $A^{*}$ is only a direct factor of the one associated with the theory $A^{*}_{\chi}$. 

Basically, all the results proved in this paper can be seen as a motivic interpretation of the computation of the rationnal $G$-theory spectrum
of a Deligne-Mumford stacks which appears in \cite{t1,t2,t3}. \\ 

In the first part of this work we will review briefly the two different definitions of Chow rings, as well as some results about the $K$-theory
of stacks which explain them. In the second and third part we will define the associated categories of motives and prove that they are
equivalent to the category of usual Chow motives. Finally, we will give some examples of motives associated to stacks. \\

\textit{Acknowledgments:} This work was motivated by a question asked by Y. Manin, who I would like to thank for pointing it out  
to me as well as for his explanations and usefull remarks. I would also like to thank A. Kresch and G. Vezzosi for their remarks.

I very grateful to Max-Planck-Institut in Bonn for hospitality and exeptional working conditions. \\

\textit{Notations:} We will work over a perfect base field $k$, of any caracteristic. An algebraic variety will be a scheme, smooth and proper
over $Spec\, k$. A $DM$-stack will be a Deligne-Mumford stack of finite type over $Spec\, k$, smooth and proper over $Spec\, k$ (\cite{dm}).
As a convention we will work in the homotopy category of stacks. Thus a morphism of stacks will be for us the class of a $1$-morphism up to 
$2$-isomorphisms. The category of $DM$-stacks thus obtained will be denoted by $\mathcal{DM}$. The full sub-category of 
varieties will be denoted by $\mathcal{VAR}$. 

\newpage

\begin{section}{Preliminaries on Chow rings of Deligne-Mumford stacks}

We start with the first definition of the Chow rings of Deligne-Mumford stacks. \\

For every $DM$-stack $F$, we consider $\mathcal{K}_{m}$, the sheaf of the $m$-th $K$-groups on $F_{et}$ associated with the
abelian presheaf
$$\begin{array}{cccc}
\mathcal{K}_{m} : & F_{et} & \longrightarrow & Ab \\
& U & \mapsto & K_{m}(U)
\end{array}$$

\begin{df}{(\cite{g})}
The codimension $m$ rationnal Chow group of $F$ is defined to be
$$A^{m}(F):=H^{m}(F_{et},\mathcal{K}_{m}\otimes \mathbb{Q}).$$
We will note $A^{*}(F):=\bigoplus_{m}A^{m}(F)$.
\end{df}

As it is shown in \cite{g}, the theory $A^{*}(F)$ is a good Chow cohomology theory. Without recalling all the properties w 
recall three of them which will be usefull for us. 

\begin{itemize}
\item (\textit{Product}) For every $DM$-stack $A^{*}(F)$ has a structure of a graded commutative ring. 

\item (\textit{Functoriality}) For every morphism of $DM$-stacks, $f : F \longrightarrow F'$, there exist an inverse image
$$f^{*} : A^{*}(F') \longrightarrow A^{*}(F)$$
which is a morphism of graded rings.

There is also a direct image
$$f_{*} : A^{*}(F) \longrightarrow A^{*}(F')$$
which is a morphism of $\mathbb{Q}$-vector spaces. This morphism is moreover graded of degree $DimF'-DimF$ if $F$ and $F'$ are pure
dimensional.

\item (\textit{Projection formula}) For every morphism $f : F \longrightarrow F'$ between two $DM$-stacks 
we have
$$f_{*}(x.f^{*}(y))=f_{*}(x).y$$
for every $x \in A^{*}(F)$ and $y \in A^{*}(F')$.

In particular, if $F$ and $F'$ are connected and of the same dimension, we have
$$f_{*}f^{*}=\times m$$
where $m$ is the generic degree of $f$ (in the stack sense \cite{v}).  

\item (\textit{Compatibility})
For every variety $X$, $A^{*}(X)\simeq CH^{*}(X)_{\mathbb{Q}}$ is the usual Chow ring of $X$.

\end{itemize}

In order to introduce the second definition $A^{*}_{\chi}$ recall the main result of \cite{t3}. It will not be used
but explains the definition of $A_{\chi}^{*}$.  \\

We begin with a $DM$-stack $F$, and consider $C_{F}^{t}$, the classifying stack of cyclic subgroups of automorphisms in $F$. 
More precisely it is defined in the following way.

A $S$-group scheme $G \longrightarrow S$ is called cyclic (and finite), if locally for the etale topology on $S$ it is isomorphic (as a
$S$-group scheme) to $\displaystyle Spec\,\frac{\mathcal{O}_{S}[T]}{T^{m}-1}$. In other words, $G$ is a multiplicative type $S$-group scheme
whose sheaf of characters has cyclic geometric fibres.

The stack $C_{F}^{t}$ is now defined by the following.

\begin{itemize}

\item For any $k$-scheme $X$, the objects in $C_{F}^{t}(X)$ are pairs $(s,c)$, where $s$ is an object in $F(X)$, and
$c$ is a sub-group scheme of the $X$-group scheme of automorphisms of $s$, $Aut_{X}(s) \longrightarrow X$, such that $c$ is
a cyclic $X$-group scheme. 

\item An isomorphism between two objects in $C_{F}^{t}(X)$, $(s,c)$ and $(s',c')$, is an isomorphism
$u : s \simeq s'$ in $F(X)$, such that $u^{-1}.c'.u=c$.

\end{itemize}

The map which sends $(s,c)$ to $s$ gives a morphism $\pi_{F} : C_{F}^{t} \longrightarrow F$, which is oubviously representable. 
Futhermore, we have the following local description of $\pi_{F}$. Locally (for the etale topology) on its moduli space 
$F$ is given by a quotient of a smooth scheme $X$ by a finite group $H$. So, to obtain a local descrition of $\pi_{F}$ it is enough
to consider the case where $F=[X/H]$.

Let $c(H)$ be a set of representative of conjugacy classes of cyclic subgroups of $H$ whose orders are prime to the caracteristic of $k$.
For each $c \in c(H)$, let $X^{c}$ the closed sub-scheme of fixed points of $c$ in $X$, and $N_{c}$ the normalisor of $c$ in $H$. 
Note that $N_{c}$ acts on $X^{c}$ by restriction. Then we have a natural equivalence
$$C_{F}^{t}\simeq \coprod_{c \in c(H)}[X^{c}/N_{c}].$$
As every $c$ is cyclic of order invertible on $X$, it is a diagonalisable group scheme, and so as $X$ is smooth, $X^{c}$ is also smooth. 

From this local description we deduce that the stack $C_{F}^{t}$ is smooth, and 
the natural morphism $\pi_{F} : C_{F}^{t} \longrightarrow F$ is representable, finite and unramified. In particular $C_{F}^{t}$ is again
smooth and proper. \\

On the stack $C_{F}^{t}$ lives the universal cyclic group stack, $q : \mathcal{C}_{F}^{t} \longrightarrow C_{F}^{t}$. It classifies
triplets $(s,c,h)$, where $(s,c)$ is an object of $C_{F}^{t}(X)$ and $h$ a section of $c$ over $X$. Thus, for any morphism
$(s,c) : U \longrightarrow C_{F}^{t}$, the pull-back of $\mathcal{C}_{F}^{t}$ on $U$ is isomorphic to the
cyclic $U$-group scheme $c \longrightarrow U$. 

Let $\chi_{F}$ be the sheaf of characters of $\mathcal{C}_{F}^{t}$ on $C_{F}^{t}$. It is defined by
$\chi_{F}:=\underline{Hom}_{Gp}(\mathcal{C}_{F}^{t},\mathbb{G}_{m})$. More explicitely, its restriction on the small etale site
of $C_{F}^{t}$ is given by
$$\begin{array}{cccc}
\chi_{F} : & (C_{F}^{t})_{et} & \longrightarrow & Ab \\
& ((s,c)\in C_{F}^{t}(U)) & \mapsto & Hom_{Gp}(c,\mathbb{G}_{m,U})
\end{array}$$
As $\mathcal{C}_{F}^{t}$ is a cyclic group stack, $\chi$ is a locally constant sheaf on $(C_{F}^{t})_{et}$, locally isomorphic
to a constant finite cyclic group sheaf.

Let us consider the sheaf of group-algebras associated to $\chi$, $\mathbb{Q}[\chi_{F}]$. It is a locally constant sheaf of $\mathbb{Q}$-algebras
on $(C_{F}^{t})_{et}$, which is locally isomorphic to the constant sheaf with fibre 
$\displaystyle \frac{\mathbb{Q}[T]}{T^{m}-1}$. As $\displaystyle \frac{\mathbb{Q}[T]}{T^{m}-1}$ is a product of cyclotomic fields with only one of maximal degree, namely
$\mathbb{Q}(\zeta_{m})$, the kernels of the local quotients 
$\displaystyle \frac{\mathbb{Q}[T]}{T^{m}-1}\longrightarrow \mathbb{Q}(\zeta_{m})$ glue together to give a well defined ideal sheaf
$\mathcal{I}_{F} \hookrightarrow \mathbb{Q}[\chi_{F}]$. We then define
$$\Lambda_{F}:=\frac{\mathbb{Q}[\chi_{F}]}{\mathcal{I}_{F}}.$$
Note that this is a well defined sheaf of $\mathbb{Q}$-algebras on $C_{F}^{t}$, locally isomorphic to the constant sheaf associated with 
a cyclotomic field. \\

We can now state the main result of \cite{t3}. For a sketch of proof the reader can consult \cite{t3}, or \cite[$3.15$]{t1} for a particular case.

\begin{thm}\label{th1}
There exist a functorial ring isomorphism
$$\phi_{F} : K_{*}(F)\otimes \mathbb{Q} \simeq H^{-*}((C_{F}^{t})_{et},\underline{K}\otimes \Lambda_{F}).$$
\end{thm}

\textit{Remark:} Here $K_{*}(F)$ is the ring of $K$-theory of perfect complexes on $F$, and $\underline{K}$ is
the presheaf of $K$-theory spectrum on $(C_{F}^{t})_{et}$. \\

The theorem justifies the following definition.

\begin{df}
For any $DM$-stack $F$, the codimension $m$ rationnal Chow group with coefficients in the characters of $F$ is defined by
$$A^{m}_{\chi}(F):=H^{m}((C_{F}^{t})_{et},\mathcal{K}_{m}\otimes \Lambda_{F}).$$
We will note $A^{*}_{\chi}(F):=\bigoplus_{m}A^{m}_{\chi}(F)$.
\end{df}

There is a natural decomposition $C_{F}^{t}\simeq F\coprod C_{F,+}^{t}$ coming from the 
section $F \longrightarrow C_{F}^{t}$ maping an object $s$ to $(s,\{e\})$. 
As the sheaf $\chi_{F}$ restricts to the constant sheaf $\mathbb{Q}$ on $F_{et}$, this induces a
group decomposition
$$A_{\chi}^{*}(F)\simeq A^{*}(F)\times A_{\chi\neq 1}^{*}(F).$$

\begin{prop}\label{p}
\begin{enumerate}

\item (Product) For every $DM$-stack, there is structure of graded commutative $\mathbb{Q}$-algebra on 
$A_{\chi}^{*}(F)$. Furthermore, the decomposition $A_{\chi}^{*}(F)\simeq A^{*}(F)\times A_{\chi\neq 1}^{*}(F)$
becomes a $\mathbb{Q}$-algebra decomposition.

\item (Functoriality) For every morphism of $DM$-stacks $f : F \longrightarrow F'$, there is an inverse image
$$f^{*} : A_{\chi}^{*}(F') \longrightarrow A_{\chi}^{*}(F)$$
which makes $A_{\chi}^{*}$ into a functor $\mathcal{DM}^{o} \longrightarrow (graded \; \mathbb{Q}-algebras)$. Furthermore 
the decomposition $A_{\chi}^{*}\simeq A^{*}\times A_{\chi\neq 1}^{*}$ is compatible with these inverse images.

There exist a direct image
$$f_{*} : A_{\chi}^{*}(F) \longrightarrow A_{\chi}^{*}(F')$$
which makes $A_{\chi}^{*}$ into a functor $\mathcal{DM} \longrightarrow \mathbb{Q}-Vect$. 

\item (Projection formula) For every morphism of $DM$-stacks $f : F \longrightarrow F'$, we have
$$f_{*}(x.f^{*}(y))=f_{*}(x).y$$
for every $x \in A_{\chi}^{*}(F)$ and $y \in A_{\chi}^{*}(F')$.

\item (Compatibility) For every variety $X$,
$A_{\chi}^{*}(X)\simeq CH^{*}(X)_{\mathbb{Q}}$ is the usual Chow ring of $X$.

\end{enumerate}
\end{prop}

\textit{Proof:} $(1)$ The product in $K$-theory gives morphisms of sheaves on $(C_{F}^{t})_{et}$
$$\mathcal{K}_{p}\otimes\mathcal{K}_{m} \longrightarrow \mathcal{K}_{p+m},$$
defining a graded ring structure on $\mathcal{K}_{*}:=\bigoplus_{m}\mathcal{K}_{m}$. By tensoring with the sheaf of
algebras $\Lambda_{F}$ we obtain a sheaf of graded $\mathbb{Q}$-algebras
$\mathcal{K}_{*}\otimes \Lambda_{F}$. It is then a general fact that the cohomology
$$A_{\chi}^{*}(F)\simeq H^{*}((C_{F}^{t})_{et},\mathcal{K}_{*}\otimes \Lambda_{F})$$
is naturally a graded $\mathbb{Q}$-algebra. \\

$(2)$ Every morphism between two $DM$-stacks $f : F \longrightarrow F'$ induces a morphism
$Cf : C_{F}^{t} \longrightarrow C_{F'}^{t}$. It sends an object $(s,c) \in C_{F}^{t}(X)$ to
$(f(s),f(c)) \in C_{F'}^{t}$. Furthermore there is a morphism of sheaves of groups
$Cf^{-1}(\chi_{F'}) \longrightarrow \chi_{F}$ given by restrictions of characters,  giving a morphism
of sheaves of algebras
$$Res_{f} : Cf^{-1}(\Lambda_{F'}) \longrightarrow \Lambda_{F}.$$
On the other hand we have inverse image in $K$-theory, which gives a morphism of sheaves on graded algebras
$$Cf^{*} : Cf^{-1}(\mathcal{K}_{*}) \longrightarrow \mathcal{K}_{*}.$$
By tensorisation this gives
$$Cf^{-1}(\mathcal{K}_{*}\otimes \Lambda_{F'}) \longrightarrow \mathcal{K}_{*}\otimes \Lambda_{F}$$
which allows to define inverse images
$$f^{*} : H^{*}((C_{F'}^{t})_{et},\mathcal{K}_{*}\otimes \Lambda_{F'}) \longrightarrow
H^{*}((C_{F}^{t})_{et},Cf^{-1}(\mathcal{K}_{*}\otimes \Lambda_{F'})) \longrightarrow  H^{*}((C_{F}^{t})_{et},\mathcal{K}_{*}\otimes \Lambda_{F}).$$

To define the direct images we use the induction morphism of characters
$$Ind_{f} : \chi_{F} \longrightarrow Cf^{-1}(\chi_{F'}).$$
This induces a morphism of sheaves of $\mathbb{Q}$-vector spaces 
$$Ind_{f} : \Lambda_{F} \longrightarrow Cf^{-1}(\Lambda_{F'}).$$
For every $m$, we use the Gersten resolution of the sheaf $\mathcal{K}_{m}$ (\cite[$7$]{g2})
$$\mathcal{K}_{m} \longrightarrow \mathcal{R}_{m}^{m} \longrightarrow \mathcal{R}_{m}^{m-1} \longrightarrow
\dots \longrightarrow \mathcal{R}_{m}^{0}.$$
Let $\mathcal{R}_{*}^{\bullet}:=\bigoplus_{m} \mathcal{R}_{m}^{\bullet}$. Thinking of $\mathcal{R}_{m}^{0}$ in cohomological
degree $0$, we have
$$A_{\chi}^{*}(F)\simeq H^{0}((C_{F}^{t})_{et},\mathcal{R}_{*}^{\bullet}\otimes \Lambda_{F}).$$
The direct image is a morphism of complexes of sheaves on $(C_{F'}^{t})_{et}$ (\cite[$7$]{g2})
$$Cf_{*}(\mathcal{R}_{*}^{\bullet}) \longrightarrow \mathcal{R}_{*}^{\bullet}.$$
Tensoring with $\Lambda_{F'}$ gives
$$Cf_{*}(\mathcal{R}_{*}^{\bullet}\otimes Cf^{-1}(\Lambda_{F'}))\simeq Cf_{*}(\mathcal{R}_{*}^{\bullet})\otimes \Lambda_{F'} \longrightarrow \mathcal{R}_{*}^{\bullet}\otimes \Lambda_{F'}.$$
We then compose with $Ind_{f}$ and take the cohomology to obtain
$$f_{*} : A_{\chi}^{*}(F) \longrightarrow A_{\chi}^{*}(F').$$

$(3)$ Using the two previous explicits definitions of $f_{*}$ and $f^{*}$ the proof is exactly the same as for
the case of scheme (\cite[$7$]{g2}). \\

$(4)$ If $X$ is a variety, then $C_{X}^{t}\simeq X$ and $\Lambda_{X}\simeq \mathbb{Q}$, so the isomorphism
$A_{\chi}^{*}(X)\simeq CH^{*}(X)_{\mathbb{Q}}$ is given by the Bloch's formula (\cite[$7$]{g2})
$$CH^{p}(X)_{\mathbb{Q}}\simeq H^{p}(X_{zar},\mathcal{K}_{p})\otimes \mathbb{Q}\simeq 
H^{p}(X_{et},\mathcal{K}_{p}\otimes \mathbb{Q}).$$ 
$\Box$ \\

\textit{Remark:} The Riemann-Roch formula of \cite[$4.11$]{t1} extends to a formula with values in $A_{\chi}^{*}$. 
Indeed, by using the construction of chern classes in \cite{g2} and the theorem \ref{th1} one can define a Chern character
$$Ch^{\chi} : K_{0}(F) \longrightarrow A_{\chi}^{*}(F).$$
The Todd class $Td(F)$ defined in \cite[$4.8$]{t1} can also be defined as $Td^{\chi}(F) \in A_{\chi}^{*}(F)$ in a very similar manner. 

To prove the Riemann-Roch formula for $Ch^{\chi}(-).Td^{\chi}$, we first use the projection formula to do galois descent and reduce
the problem to the case where $k$ is algebraically closed. We choose an embending $\mu_{\infty}(k) \hookrightarrow
\mathbb{C}^{*}$. 
Then the formula follows from \cite[$3.36$]{t2} and the fact that (see \ref{l} for
a proof of this)
$$A_{\chi}^{*}(F)\otimes \mathbb{Q}(\mu_{\infty}(k)) \simeq A_{rep}^{*}(F).$$
It is also true that the Chern character
$$Ch^{\chi} : K_{0}(F)_{\mathbb{Q}} \longrightarrow A^{*}_{\chi}(F)$$
is a ring isomorphism.

\end{section}

\begin{section}{First construction}

The construction of the category of Chow motives for $DM$-stacks using the theory $A^{*}$ was done in \cite[$8$]{bm}. We will 
denote it by $\mathcal{M}^{DM}$, and call its objects the $DMC$-motives, as suggested in \cite{bm}. We start by recalling briefly its construction. \\

For $F, F' \in\mathcal{DM}$, we define the vector space of correspondences of degree $m$ between $F$ and $F'$
$$S^{m}(F,F'):=\{x \in A^{*}(F\times F') / (p_{2})_{*}(x) \in A^{m}(F')\}$$
where $p_{2} : F\times F' \longrightarrow F'$ is the second projection. 

We have the usual composition
$$\circ : S^{m}(F,F')\otimes S^{n}(F',F'') \longrightarrow S^{p+m}(F,F'')$$
given by the formula 
$$x\circ y:=(p_{13})_{*}(p_{12}^{*}(x).p_{23}^{*}(y)),$$
where the $p_{ij}$ are the natural projections of $F\times F'\times F''$ on two of the three factors. \\

Objects of $\mathcal{M}^{DM}$ are triplets $(F,p,m)$, with $F \in \mathcal{DM}$, $p$ an idempotent in the ring
of correspondences $S^{0}(F,F)$, and $m \in \mathbb{Z}$. 
The morphisms between $(F,p,m)$ and $(F',q,n)$ are given by
$$Hom_{\mathcal{M}^{DM}}((F,p,m),(F',q,n)):=q\circ S^{n-m}(F,F')\circ p \subset S^{n-m}(F,F').$$

Recall also that for any morphism in $\mathcal{DM}$, $f : F \longrightarrow F'$, we have its graph
$$\Gamma_{f}=f\times Id : F \longrightarrow F'\times F,$$
and so a well defined element
$$[f^{*}]:=(\Gamma_{f})_{*}(1) \in S^{0}(F',F).$$
We can also consider its transposed 
$$[f_{*}]:=[f^{*}]^{t} \in S^{*}(F,F').$$
This allows us to define a functor
$$\begin{array}{cccc}
h : & \mathcal{DM}^{o} & \longrightarrow & \mathcal{M}^{DM} \\
& F & \mapsto & (F,Id,0) \\
& f & \mapsto & [f^{*}]
\end{array}$$

As for the case of varieties, the category $\mathcal{M}^{DM}$ is $\mathbb{Q}$-linear karoubian category (\cite[$8.1$]{bm}). In particular, 
this implies that if a morphism in $\mathcal{M}^{DM}$ possesses a left inverse then it is a direct factor.

It is also symetric monoidal for the tensor product defined by
$$(F,p,m)\otimes (F',q,n):=(F\times F',p\otimes q,n+m).$$

As usual we shall write $\mathbb{L}^{m}=(Spec\, k,Id,m)$ for the $m$-th power of the Lefschetz motive. 
Note that for every $DMC$-motive $M$, we have $M\simeq (F,p,0)\otimes \mathbb{L}^{m}$. As $(F,p,0)$ is a direct factor
in $h(F)$, this shows that $M$ is a direct factor in 
some $h(F)\otimes \mathbb{L}^{m}$. 

Finally, there is a natural fully faithfull tensorial functor
$$\mathcal{M} \longrightarrow \mathcal{M}^{DM}$$
from the usual category of Chow motives of varieties to the category of $DMC$-motives. This functor fits into a commutative
diagramm
$$\xymatrix{
\mathcal{VAR}^{o} \ar[r] \ar[d]_-{h} & \mathcal{DM}^{o} \ar[d]^-{h} \\
\mathcal{M} \ar[r] & \mathcal{M}^{DM} }$$

The following theorem is a positive answer to the question \cite[$8.2$]{bm}.

\begin{thm}\label{th2}
The previous functor
$$\mathcal{M} \longrightarrow \mathcal{M}^{DM}$$
is an equivalence of $\mathbb{Q}$-tensorial categories.
\end{thm}

\textit{Proof:} By noticing that the essential image is closed by direct factors 
(because any direct factor of $(X,p,m)$ in $\mathcal{M}^{DM}$ is of the form $(X,p \circ q,m)$), we only have to check that
for each connected $F \in \mathcal{DM}$, $h(F)$ is a direct factor of some $h(X)$ for $X \in \mathcal{VAR}$. \\

Let $F \in \mathcal{DM}$, and by \cite[$16.6$]{e} choose an integral scheme $X$ and a finite and surjective morphism
$X \longrightarrow F$. Using \cite{dj} we can find $Y \longrightarrow X$ which is generically finite, with $Y$ a variety.
We know consider the composed morphism $f : Y \longrightarrow F$, as well as
$$[f_{*}] : h(Y) \longrightarrow h(F)$$
$$[f^{*}] : h(F) \longrightarrow h(Y).$$
The indentity principle \cite[$8.2$]{bm} and the projection formula implies that $\frac{1}{m}.[f_{*}]$ is a left inverse to 
$[f^{*}]$. This implies that $[f^{*}]$ is a direct factor. More explicitely we have
$h(F)\simeq (X,\frac{1}{m}[f_{*}]\circ [f^{*}],0)$. $\Box$ \\

Inverting the equivalence $\mathcal{M} \longrightarrow \mathcal{M}^{DM}$ gives a functor
$$h : \mathcal{DM}^{o} \longrightarrow \mathcal{M}.$$
As an inverse of a monoidal functor has a natural monoidal structure, $h$ is naturally a monoidal functor. We obtain this way 
natural isomorphisms
$$h(F\times F')\simeq h(F)\otimes h(F')$$
wich are associatives, commutatives and unitaries. In particular, the diagonal of a $DM$-stack $F$ gives a commutative algebra structure
on the motive $h(F)$. This can be used for example to show that every good cohomology theory for 
varieties extends to $DM$-stacks.

\end{section}

\begin{section}{Second construction}

\begin{df}
For two $DM$-stacks $F$ and $F'$, we define the vector space of $\chi$-correspondences of degree $m$ between $F$ ans $F'$
by
$$S^{m}_{\chi}(F,F'):=\{x \in A_{\chi}^{*}(F\times F') / (p_{2})_{*}(x) \in A_{\chi}^{m}(F')\}.$$
\end{df}

As in the previous case, we have a composition
$$\circ : S_{\chi}^{m}(F,F')\otimes S_{\chi}^{n}(F',F'') \longrightarrow S_{\chi}^{p+m}(F,F'')$$
given by the formula
$$x\circ y:=(p_{13})_{*}(p_{12}^{*}(x).p_{23}^{*}(y)),$$
where the $p_{ij}$ are the natural projections of $F\times F'\times F''$ on two of the three factors. 

\begin{df}
We define the category of $DMC_{\chi}$-motives, $\mathcal{M}^{DM}_{\chi}$ as follows.

\begin{itemize}
\item Objects of $\mathcal{M}^{DM}_{\chi}$ are triplets $(F,p,m)$, where $F$ is a $DM$-stack, $p \in S^{0}_{\chi}(F,F)$ an
idempotent, and $m \in \mathbb{Z}$. 

\item The set of morphisms between $(F,p,m)$ and $(F',q,n)$ is defined by
$$Hom_{\mathcal{M}^{DM}_{\chi}}((F,p,m),(F',q,n)):=
q\circ S^{n-m}_{\chi}(F,F')\circ p \subset S^{n-m}_{\chi}(F,F').$$

\item The composition of morphisms in $\mathcal{M}^{DM}_{\chi}$ is given by composition of
$\chi$-correspondences.
\end{itemize}
\end{df}

For any morphism $f : F \longrightarrow F'$ between two $DM$-stacks, we define  as usual
$$[f^{*}]:=(\Gamma_{f})_{*}(1) \in S^{0}_{\chi}(F',F),$$
as well as its transposed
$$[f_{*}]:=[f^{*}]^{t} \in S^{*}_{\chi}(F,F').$$
Using this we define a natural functor
$$\begin{array}{cccc}
h_{\chi} : & \mathcal{DM}^{o} & \longrightarrow & \mathcal{M}_{\chi}^{DM} \\
& F & \mapsto & (F,Id,0) \\
& f & \mapsto & [f^{*}]
\end{array}$$

The same arguments as for motives of varieties show that $\mathcal{M}^{DM}_{\chi}$ is a $\mathbb{Q}$-tensorial karoubian 
category. There is alos a tensor product, given as usual by $(F,p,m)\otimes (F',q,n):=(F\times F',p\otimes q,m+n)$. 

Note that the compatibility propety of \ref{p} implies that there is a natural fully faithfull functor
$$\mathcal{M} \longrightarrow \mathcal{M}^{DM}_{\chi}.$$

\begin{df}
For any $DMC_{\chi}$-motive $M$, we define its $m$-th Chow group by
$$A^{m}_{\chi}(M):=Hom_{\mathcal{M}^{DM}_{\chi}}(\mathbb{L}^{m},M).$$
We will note $A^{*}_{\chi}(M):=\bigoplus_{m}A^{m}_{\chi}(M)$.
\end{df}

\textit{Remark:} Using the Chern character we have
$$Ch^{\chi} : K_{0}(F)_{\mathbb{Q}}\simeq A_{\chi}^{*}(F).$$

Finally, the indentity principle says that the functor
$\mathcal{M}^{DM}_{\chi} \longrightarrow Hom(\mathcal{DM}^{o},Ab)$, which sends $M$ to the functor $F \mapsto A^{*}_{\chi}(M\otimes h(F))$ 
is fully faithfull (it follows immediately from the Yoneda lemma and the fact that every $DMC_{\chi}$-motive
is a direct factor of a $h(F)\otimes \mathbb{L}^{m}$).

\begin{thm}\label{th3}
The natural functor
$$\mathcal{M} \longrightarrow \mathcal{M}^{DM}_{\chi}$$
is an equivalence of $\mathbb{Q}$-tensorial categories.
\end{thm}

\textit{Proof:} As for \ref{th2} it is enough to show that every $h_{\chi}(F)$ is a direct factor in some $h_{\chi}(X)$. \\

For the next lemma recall that for any stack $F$ we can define its inertia stack $I_{F}$ (\cite{v}), whose objects are pairs
$(s,h)$, with $s$ an object of $F$ and $h$ and automorphism of $s$. It can for example be defined by the formula
$$I_{F}:=F\times_{F\times F}F.$$

\begin{lem}\label{l0}
For any Deligne-mumford stack proper over $Spec\, k$ (non necerally smooth), there exist varieties $Y_{i}$ and 
finite groups $H_{i}$ together with a proper representable morphism
$$F_{0}:=\coprod_{i}Y_{i}\times BH_{i} \longrightarrow F$$
such that the induced morphism
$$I_{F_{0}} \longrightarrow I_{F}$$
is generically finite and surjective.
\end{lem}

\textit{Proof:} By \cite[$16.6$]{e} 
we can choose a finite and surjective morphism $X \longrightarrow F$, with $X$ a normal scheme. Let $F \longrightarrow M$
be the moduli space of $F$, and consider $F_{X}$, the normalization of the fibre product $F\times_{M}X$. By definition, the stack $F_{X}$
is normal and the projection to its moduli space $F_{X} \longrightarrow X$ possesses a section. It follows from \cite[$2.7$]{v}
that $F_{X}$ is a neutral gerb. By choosing a finite and etale morphism $Y \longrightarrow X$ and defining $F':=F_{X}\times_{X}Y$, we 
find a trivial gerb $F'\simeq Y\times BH$, together with a morphism $F' \longrightarrow F$. By construction this morphism is generically obtained
by a pull back of a etale morphism on $M$. This implies that there exists a dense
open sub-stack $U$ of $F$, such that $I_{U} \subset I_{F}$ is in the image of $I_{F'} \longrightarrow I_{F}$. Proceding by noetherian induction
we find reduced schemes $X_{i}$, and finite groups $H_{i}$, with a morphism $F':=\coprod_{i}X_{i}\times BH_{i} \longrightarrow F$ such 
that $I_{F'} \longrightarrow I_{F}$ is finite and surjective. 

We now apply \cite{dj} to each $X_{i}$ and choose generically finite morphism $Y_{i} \longrightarrow X_{i}$, with $Y_{i}$ a variety. Let
$F_{0}:=\coprod_{i}Y_{i}\times BH_{i}$. As $I_{F_{0}}\simeq I_{F}\times_{F}F_{0}$, the induced morphism
$$I_{F_{0}} \longrightarrow I_{F}$$
is still surjective and generically finite. $\Box$ \\

\begin{lem}\label{l}
Suppose that $k$ is algebraically closed, and choose an embedding $\mu_{\infty}(k) \hookrightarrow \mathbb{C}^{*}$.
Let $F$ be a $DM$-stack, and denote by $I_{F}^{t}$ the open and closed sub-stack of $I_{F}$ 
whose objects are pairs $(s,h)$, such that the order of $h$ is prime to the characteristic of $k$.
Then there exist an $\mathbb{Q}(\mu_{\infty}(k))$-algebra isomorphism
$$A^{*}_{\chi}(F)\otimes \mathbb{Q}(\mu_{\infty}(k))\simeq A^{*}(I_{F}^{t})\otimes \mathbb{Q}(\mu_{\infty}(k)).$$
Furthermore this isomorphism is compatible with inverse and direct images.
\end{lem}

\textit{Proof:} Let $u : I_{F}^{t} \longrightarrow C_{F}^{t}$ the morphism which sends an object $(s,h)$ to $(s,<h>)$, where
$<h>$ is the subgroup generated by $h$ in $Aut(s)$. This is a representable finite et etale morphism.
It is easy to see that there is an isomorphism of 
sheaves of graded $\mathbb{Q}(\mu_{\infty}(k))$-algebras on $(C_{F}^{t})_{et}$
$$u_{*}(\mathcal{K}_{*}\otimes \mathbb{Q}(\mu_{\infty}(k)))\simeq \mathcal{K}_{*}\otimes \Lambda_{F}\otimes \mathbb{Q}(\mu_{\infty}(k)).$$
This induces the required isomorphism
$$A^{*}(I_{F}^{t})\otimes \mathbb{Q}(\mu_{\infty}(k)) \simeq
H^{*}((I_{F}^{t})_{et},\mathcal{K}_{*}\otimes \mathbb{Q}(\mu_{\infty}(k))) \simeq
H^{*}((C_{F}^{t})_{et},u_{*}(\mathcal{K}_{*}\otimes \mathbb{Q}(\mu_{\infty}(k))))$$
$$\simeq H^{*}((C_{F}^{t})_{et},\mathcal{K}_{*}\otimes \Lambda_{F}\otimes \mathbb{Q}(\mu_{\infty}(k)))\simeq 
A^{*}_{\chi}(F)\otimes\mathbb{Q}(\mu_{\infty}(k)).$$
The compatibility with inverse and direct images is clear by definitions. $\Box$ \\

Let $g : F_{0}:=\coprod_{i}Y_{i}\times BH_{i} \longrightarrow F$ be a morphism as in \ref{l0}.

\begin{lem}
The element $\beta:=g_{*}(1) \in A^{*}_{\chi}(F)$ is invertible. 
\end{lem}

\textit{Proof:} We first use the projection formula \ref{p} to show that for any finite extension $k'/k$,
we have a natural isomorphism of algebras
$$A^{*}(F)_{\chi}\simeq A^{*}_{\chi}(F\times_{Spec\, k} Spec\, k')^{Gal(k'/k)}.$$
This allows to assume that $k$ is algebraically closed. 

Applying the lemma \ref{l}, it is enough to show that $Ig_{*}(1) \in A^{*}(I_{F}^{t})$ is invertible. 
But, as $Ig$ is generically finite and surjective this is obvious. $\Box$ \\

Consider $\Delta_{*}(\beta) \in A^{*}_{\chi}(F\times F)=S^{*}_{\chi}(F,F)$, where $\Delta : F \longrightarrow
F\times F$ is the diagonal. By the previous lemma $\beta$ is invertible in the graded ring
$S^{*}(F,F)$. Let $\alpha:=[g_{*}]\circ \beta^{-1} \in S^{*}(F',F)$. Then we have
$[g^{*}]\circ \alpha=1$. This shows that the $0$-th component of $\alpha$ is a left inverse
of $[g^{*}]$. As the category $\mathcal{M}^{DM}_{\chi}$ is karoubian, this implies that $[g^{*}]$ is a direct factor, and so that $h_{\chi}(F)$ is 
a direct factor in $h_{\chi}(F')$. \\

As $h_{\chi}(F')\simeq \bigoplus_{i}h_{\chi}(X_{i})\otimes h_{\chi}(BH_{i})$ it remains to show that for any finite group $H$, 
$h_{\chi}(BH)$ is isomorphic to some power of the trivial motive $h_{\chi}(Spec\, k)$. 

Let $Ch^{\chi} : K_{0}(BH) \longrightarrow A^{0}_{\chi}(BH)$ the Chern character, $\rho_{1}, \dots, \rho_{r}$ a set of representatives
of irreducibles representations of $H$ over $k$, and $\alpha_{i}:=Ch^{\chi}(\alpha_{i})$. These elements define morphisms
of $DMC_{\chi}$-motives
$\alpha_{i} : h_{\chi}(Spec\, k) \longrightarrow h_{\chi}(BH)$
Let us consider the sum 
$$\bigoplus_{i}\alpha_{i} : h_{\chi}(Spec\, k)^{r} \longrightarrow h_{\chi}(BH),$$
and prove that it is an isomorphism.
By the identity principle, we have to show that for every $DM$-stack $F$, the induced morphism
$$\bigoplus_{i}\alpha_{i} : (A^{*}_{\chi}(F))^{r} \longrightarrow A^{*}_{\chi}(F\times BH)$$
is an isomorphism. But as $Ch^{\chi}$ is an isomorphism, the previous morphism is isomorphic to the Kunneth morphism
$$A^{*}_{\chi}(F)\otimes A^{*}_{\chi}(BH) \longrightarrow A^{*}_{\chi}(F\times BH),$$
and so the theorem follows from the following lemma.

\begin{lem}
For every $DM$-stack $F$ and every finite group $H$, the Kunneth morphism
$$A^{*}_{\chi}(F)\otimes A^{*}_{\chi}(BH) \longrightarrow A^{*}_{\chi}(F\times BH)$$
is an isomorphism.
\end{lem}

\textit{Proof:} Using galois descent we can suppose that $k$ is algebraically closed. Then, using the lemma \ref{l} we reduce the problem to show
that the Kunneth morphism
$$A^{*}(I_{F}^{t})\otimes A^{*}(I_{BH}^{t}) \longrightarrow A^{*}(I_{F\times BH}^{t})$$
is an isomorphism. 

Let $A$ be a set of representative of conjugacy classes of elements in $H$ with order prime to the characteristic of $k$. We have
$$I_{BH}^{t}\simeq \coprod_{h \in A}BZ_{h} \qquad I_{F\times BH}^{t}\simeq I_{F}^{t}\times I_{BH}^{t},$$
where $Z_{h}$ is the centralisator of $h$ in $H$. So we only need to prove that the Kunneth morphism
$$A^{*}(I_{F}^{t})\otimes A^{*}(BZ_{h}) \longrightarrow A^{*}(I_{F}^{t}\times BZ_{h})$$
is an isomorphism.
But this morphism fits into a commutative diagram
$$\xymatrix{
A^{*}(I_{F}^{t})\otimes A^{*}(BZ_{h}) \ar[r] & A^{*}(I_{F}^{t}\times BZ_{h}) \\
A^{*}(I_{F}^{t}) \ar[u]^-{Id \otimes 1} \ar[ru]_-{v^{*}} & }$$
where $v : I_{F}^{t}\times BZ_{h} \longrightarrow I_{F}^{t}$ is the first projection.  
Now, as $A^{*}(BZ_{h})\simeq \mathbb{Q}$, the vertical morphism is an isomorphism. On the other hand, $v$ has a natural section
$u : I_{F}^{t} \longrightarrow I_{F}^{t}\times BZ_{h}$ and the projection formula  
shows that $u^{*}$ is an isomorphism, which implies that $v^{*}$ is an isomorphism. $\Box$ \\

Inverting the equivalence $\mathcal{M} \longrightarrow \mathcal{M}^{DM}_{\chi}$ gives a functor
$$h_{\chi} : \mathcal{DM}^{o} \longrightarrow \mathcal{M}.$$
As for the case of the first construction this functor has a natural monoidal structure. This implies that for any 
$DM$-stack $F$, the motive $h_{\chi}(F)$ has a natural structure of a commutative algebra in $\mathcal{M}$. In particular any
good cohomology theory for varieties extends trough $h_{\chi}$ to a new theory for stacks. 

\begin{prop}
The functor $h$ is a direct factor of the functor $h_{\chi}$.
\end{prop}

\textit{Proof:} This follows immediately from the natural decomposition 
$A^{*}_{\chi}\simeq A^{*}\times A_{\chi\neq 1}^{*}$. $\Box$ \\

\textit{Remark:} For any complex variety $V$ and $\beta \in H_{2}(V,\mathbb{Z})$, we can define the Gromov-Witten
correspondence (\cite{b})
$$I_{g,n}(V,\beta) \in S^{*}(V^{n},\mathcal{M}_{g,n}),$$
which is a morphism of graded $DMC$-motives (\cite[$8$]{bm}). It seems natural to ask if this correspondence extends in a natural way to
$$I_{g,n}^{\chi}(V,\beta) \in S^{*}_{\chi}(V^{n},\mathcal{M}_{g,n})$$
(i.e. as a morphism of graded $DMC_{\chi}$-motives). This question is of course linked to the question of constructing an extended 
virtual fundamental class $\mathcal{I}_{g,n}^{\chi}(V,\beta) \in A_{*}^{\chi}(\mathcal{M}_{g,n}(V,\beta))$. 

\end{section}

\begin{section}{Examples}

We have seen that the two Chow coholomogy theories $A^{*}$ and $A^{*}_{\chi}$ give natural functors
$$h, h_{\chi} : \mathcal{DM}^{o} \longrightarrow \mathcal{M},$$
such that $h$ is a direct factor of $h_{\chi}$. In this last chapter we will give some examples of motives associated to certain stacks, and see
some expicit relations between $h_{\chi}$ and $h$. 

The proofs of the following three facts are left to the reader (they all follow from the indentity principle and
the explicit description of the stacks $C_{F}^{t}$ and the sheaves $\Lambda_{F}$). For the sake of simplicity we will suppose
that $k$ contains the roots of unity. \\

If a finite group $H$ acts on a motive $M$ we will denote by $M^{H}$ the direct factor of $M$
corresponding to the projector $\displaystyle \frac{1}{m}.\sum_{h \in H}h$. 

\begin{enumerate}

\item \textit{Quotients stacks} \\

Let $H$ be a finite group acting on a variety $X$. Let $c(H)$ be a set of representatives of conjugacy classes of cyclic
sub-groups of $H$, whose orders are prime to the characteristic of $k$. For every $c \in c(H)$ let $X^{c}$ be the sub-variety of $X$ of
fixed points of $c$, and $N_{c}$ the normalisator of $c$ in $H$. For any $c \in c(H)$, let 
$s(c)$ be the set of injectives characters $c \longrightarrow k^{*}$.

Then the group $N_{c}$ acts on $X^{c}$ and $s(c)$, and so on the product $X^{c}\times s(c)$,
and there is an isomorphism
$$h_{\chi}([X/H])\simeq \bigoplus_{c \in c(H)}h(X^{c}\times s(c))^{N_{c}}.$$
Furthermore, the motive $h([X/H])$ corresponds to the component of the trivial sub-group
$$h([X/H])\simeq h(X)^{H}.$$ 

For example if $X=Spec\, k$ we obtain
$$h_{\chi}(BH)\simeq h(Spec\, k)^{r},$$
where $r$ is the number of irreducible representations of $H$ in $k$-vetor spaces. But notice that this isomorphism does not preserve
the product structures (given on any $h_{\chi}(F)$ by the diagonal morphism). Indeed, if $\rho_{1}, \dots, \rho_{r}$ are the
irreducibles representations of $H$ over $k$, then we have the mutiplication rules
$$\rho_{i}\otimes \rho_{j}\simeq \bigoplus_{k}\rho_{k}^{n^{i,j}_{k}}.$$
Then, the product on $h_{\chi}(BH)$ corresponds on $h(Spec\, k)^{r}$ to the morphism 
$$h(Spec\, k)^{r}\otimes h(Spec\, k)^{r}\simeq h(Spec\, k)^{r^{2}} \longrightarrow h(Spec\, k)^{r}$$
given by the $r^{2}$ by $r$ matrix $(n^{i,j}_{k})_{i,j,k}$. 

\item \textit{Gerbs} \\

Let $F$ be a connected $DM$-stack which is a gerb (i.e. the morphism $C_{F}^{t} \longrightarrow F$ is etale), and $F \longrightarrow X$ its 
projection to its moduli
space. Recall that locally for the etale topology of $X$, $F$ is equivalent to $X\times BH$, for $H$ a finite group. This defines a locally constant
sheaf of groups up to inner automorphisms on $X_{et}$, which is classified by its monodromy
$$\pi_{1}^{et}(X) \longrightarrow Out(H).$$
Let $cycl(H)$ be the set of cyclic sub-groups of $H$ of order prime to the characteristic of $k$, and for any $c \in cycl(H)$, $s(c)$ the
set the of injectives characters $c \longrightarrow k^{*}$. The group $H$ acts by conjugaison on
$\displaystyle \coprod_{c \in cycl(H)}s(c)$, and let
$\displaystyle R(H):=(\coprod_{c \in cycl(H)}s(c))/H$ be the quotient. 
The group $Aut(H)$ acts naturally on $R(H)$ and
any inner automorphisms of $H$ acts trivially, so we deduce a morphism
$$\pi_{1}^{et}(X) \longrightarrow Aut(R(H)),$$
which it turns corresponds to a finite etale covering $Y \longrightarrow X$. 

There is then an isomorphism
$$h_{\chi}(F)\simeq h(Y).$$
Note that the trivial subgroup with the trivial character induces a section $X \longrightarrow Y$, which gives a decomposition
$$h(Y)\simeq h(X)\oplus h(Y)^{\neq 1}.$$
Furthermore, $h(F)$ corresponds to the factor $h(X)$.

\item \textit{$1$-Dimensional complex orbifolds} \\

Suppose that $k=\mathbb{C}$ is the field of complex number, and that $F$ is a $1$-dimensional $DM$-stack, which is generically a variety (i.e.
$C_{F}^{t} \longrightarrow F$ is birationnal). Let $C$ be the moduli space of $F$, which is a smooth projective curve, and note
$x_{1}, \dots, x_{r}$ the points of $C$ where $F$ is not a scheme. Locally for the analytic topology around each $x_{i}$, $F$ is a quotient stack of a disc by a cyclic group $\mathbb{Z}/n_{i}$. 
There is then an isomorphism
$$h_{\chi}(F)\simeq h(C)\bigoplus_{i}h(Spec\, \mathbb{C})^{n_{i}-1},$$
where $h(F)$ corresponds to the factor $h(C)$. 

\end{enumerate}

\end{section}

\end{document}